\documentclass[twoside,11pt]{article}
\usepackage{graphicx}
\usepackage{amssymb}
\usepackage{amsthm}
\newfont{\bbb}{msbm10 scaled\magstep 1}
\newcommand{\bd}{{\rm bd}}
\newcommand{\inter}{{\rm int}}
\newcommand{\diam}{{\rm diam}}

\newcommand{\width}{{\rm width}}
\newtheorem{thm}{Theorem}
\newtheorem{cor}{Corollary}
\newtheorem{pro}{Proposition} 
\newtheorem{cla}{Claim}

\font\bigbold=cmbx10 at 14 pt
\font\bbigbold=cmbx10 at 16 pt 
\date{}

\title
{\bbigbold 
Application of spherical convex bodies 
to Wulff shape}

\begin{document}

\baselineskip 17.25 pt

\maketitle

\vskip -1.3cm 
\centerline
{\bigbold Marek Lassak}

\vskip0.25cm
\centerline {\it University of Science and Technology}

\centerline {\it Kaliskiego 7, 85-789 Bydgoszcz, Poland}

\centerline {\it e-mail: lassak@utp.edu.pl}

\vskip0.65cm

\pagestyle{myheadings} \markboth{ Marek Lassak}{Application of spherical convex bodies to Wulff shape}

\vskip 0.4cm

\noindent
{\bf Abstract.}
After a few claims on lunes and spherical convex sets we present some relationships between the diameter, width and thickness of reduced spherical convex bodies and bodies of constant diameter.
These relationships are formulated and proved in order to apply them for the final theorem, which permits to recognize if a Wulff shape in the Euclidean space is self-dual.



\date{}

\maketitle

\vskip 0.2cm
\section{Introduction}

Our subject is the spherical geometry (see the monographs \cite{Pa} and \cite{VB}) and also a theorem on recognizing if a Wulff shape in the Euclidean space is self-dual.
We start with necessary notions and proving a number of claims, theorems and propositions concerning spherical geometry which, as also some results from \cite{L-AEQ}, \cite {LaMu-BULL}, \cite{LaMu-AEQ} and \cite {LaMu-FAS}, are needed for this  theorem on Wulff shape.

In $E^{d+1}$, where $d\geq 2$, take the unit sphere $S^d$ centered at the origin. 
The intersection of $S^d$ with any $(k+1)$-dimensional Euclidean subspace, where $0 \leq k \leq d-1$, is called a {\it $k$-dimensional subsphere of $S^d$}.
For $k=1$ we call it a {\it great circle}, and for $k=0$ a {\it pair of antipodes}.
If different points $a, b \in S^d$ are not antipodes, by the {\it arc} $ab$ connecting them we mean this part of the great circle containing $a$ and $b$ which does not contain any pair of antipodes. 
By the {\it spherical distance} $|ab|$, or shortly {\it distance}, of these points we understand the length of the arc connecting them. 

The intersection of $S^d$ with any half-space of $E^{d+1}$ is called a {\it hemisphere} of $S^d$.
In other words, a {\it hemisphere} $H(c)$ of $S^d$ is the set of points of $S^d$ in distances at most $\frac{\pi}{2}$ from a point $c$ called the {\it center} of this hemisphere.
Two hemispheres whose centers are antipodes are called {\it opposite hemispheres}.

We say that a set $C \subset S^d$ not containing any pair of antipodes is {\it convex} if together with every two its points it contains the arc connecting them.
If the interior $\inter (C)$ of a closed convex set $C \subset S^d$ is non-empty, $C$ is called a {\it convex body}.
We call $C$ {\it strictly convex} if in its boundary $\bd (C)$ there is no arc.
If  $A$ is a subset of a convex set of $S^d$, then by the {\it convex hull of $A$} we mean the intersection of all convex sets containing $A$ (so it is the smallest convex superset of~$A$).
We call $e$ an extreme point of a convex body $C \subset S^d$ provided $C \setminus \{e\}$ is convex.

If a hemisphere $H$ contains a convex body $C$ and if $p \in \bd (H) \cap C$, we say that $H$ {\it supports $C$ at $p$}.
We also say that $H$ is a {\it supporting hemisphere of $C$ at $p$}.
If at every boundary point of a convex body $C \subset S^d$ exactly one hemisphere supports $C$, then $C$ is called {\it smooth}.

If hemispheres $G$ and $H$ of $S^d$ are different and not opposite, then $L = G \cap H$ is called {\it a lune} of $S^d$. 
This notion is considered in many books and papers.  
The parts of $\bd (G)$ and $\bd (H)$ contained in $G\cap H$ are denoted by $G/H$ and $H/G$, respectively.
Clearly, $(G/H) \cup (H/G)$ is the boundary of the lune $G \cap H$.
Points of $(G/H) \cap (H/G)$ are called {\it corners} of the lune $G \cap H$. 
The set of them is denoted by ${\rm corn} (L)$.
By the {\it thickness $\Delta (L)$ of the lune} $L = G \cap H$ we mean the spherical distance of the centers of $G/H$ and $H/G$.

Section 2 gives claims on lunes and spherical convex sets.
Sections 2 and 3 recall the spherical notions of width, bodies of constant width and constant diameter, and
reduced bodies. 
Section 3 presents some new relationships between the thickness and the diameter of reduced bodies.
Here we also show that every reduced body of thickness at least $\frac{\pi}{2}$ is of constant width.
Section 4 is devoted to applications of these facts for recognizing if a Wulff shape is self-dual.

\section{Four claims on lunes and spherical convex sets}

The following claim is obvious. 

\begin{cla}  \label{twolunes}
Let $H(c)$ be any hemisphere of $S^d$.
Then any $(d-1)$-dimensional subsphere of $S^d$ containing $c$ dissects $H(c)$ into two lunes of thickness $\frac{\pi}{2}$.
\end{cla}

For any convex body $C \subset S^d$ and any hemisphere $K$ supporting $C$ we define the {\it width of $C$ determined by $K$} as the minimum thickness of a lune $K \cap K'$ over all hemispheres $K' \not = K$ supporting $C$ and we denote it by ${\rm width}_K (C)$. 
By the {\it thickness} $\Delta (C)$ of $C$ we mean the minimum  of ${\rm width}_K (C)$ over all hemispheres $K$ supporting $C$.
In the literature $\Delta (C)$ is also called {\it the minimum width of $C$}.
Clearly, $\Delta (C)$ is nothing else but the thickness of each ``narrowest'' lune containing $C$. 
We say that $C$ is of {\it constant width} $w$ provided all its widths ${\rm width}_K (C)$ are equal to $w$.
These notions and a few properties of lunes and convex bodies on $S^d$ are presented in \cite{L-AEQ} and \cite{LaMu-AEQ}.
The book \cite{MMO} gives a wide survey of results on bodies of constant width in various structures.
Here is an additional property needed later. 

\begin{cla} \label{diameter}
For every convex body $C \subset S^d$ we have $\Delta (C) \le \diam (C)$. 
Moreover, if $\Delta(C) = \diam (C)$ and if they are most $\frac{\pi}{2}$, then $C$ is of constant width $w$ equal to $\Delta(C)= \diam (C)$.
\end{cla}

\begin{proof}
We get the first assertion by the definition of $\Delta (C)$ and the inequality 
$$\max \{\width_K(C): K \ {\rm is \ a \ supporting \ hemisphere \ of} \ K  \} \le \diam (C)$$ 
resulting from Theorem 3 and Proposition 1 of \cite{L-AEQ}.

In order to show the second assertion, let us apply Theorem 3 of \cite{L-AEQ} which says that if $\diam (C)$ is at most $\frac{\pi}{2}$, then it is nothing else but the maximum of the widths of $C$. 
Moreover, having in mind that $\Delta (C)$ is the minimum width of $C$, from the assumption that $\Delta (C) =\diam (C) = w$ we conclude that all widths of $C$ are equal to $w$, which means that $C$ is of constant width~$w$.
\end{proof}

\begin{cla}  \label{corners}
Let $L \subset S^d$ be a lune. 
Let $C \subset L$ be a convex body such that the set $F= C \cap {\rm corn} (L)$ is non-empty.  
Then at least one extreme point of $C$ is in ${\rm corn} (L)$.
\end{cla}

\begin{proof}
Take the smallest subsphere $S^k$ of $S^d$ containing $F$.
Of course, $k \geq 1$ and $F$ is a convex body on $S^k$.
Thus $F$ is the convex hull of its extreme points (see \cite{L-AEQ}, p. 565).

Clearly, the convex body $F$ of $S^k$ has at least one extreme point $e$. 
So for every $p, q \in F$ with $e \in pq$, our point $e$ must be an end-point of $pq$.

Our aim is to show that $e$ is an extreme point of $C$, as well.
So it is sufficient to show that always if $e \in ab$ for $a, b \in C$, then $e$ is an end-point of $ab$.
Thus assume that $a,b \in C$ and $e \in ab$.

Case 1, when $a, b \in F$.
By the second paragraph of this proof, $e$ is an end-point of $ab$. 

Case 2, when $a \in F$ and $b \in C \setminus F$, or vice-versa.
Since 
$C \setminus F$ is a convex set, all points of the arc $ab$, besides $a$, are in 
$C \setminus F$ and $a \in F$.
Also recall that $e \in F$.
Thus if $e \in ab$, then $e=a$. 

Case 3, when $a, b \in C setminus F$. 
Then since $C \setminus F$ is convex, we see that
 $ab \subset C \setminus F$. 
So our $e$ (which belongs to $F$) cannot be in $ab$,which means that this case is impossible.

From these cases we conclude that $e$ is an extreme point of $C$.
\end{proof}

\begin{cla} \label{extreme}
Let $C \subset S^d$ be a convex body with. 
If $\diam (C) = \frac{\pi}{2}$, then there are two points of $C$ in the distance $\frac{\pi}{2}$ such that at least one of them is an extreme point of $C$. 
If $\diam (C) < \frac{\pi}{2}$, then every two points of $C$ in the distance $\diam (C)$ are extreme.
\end{cla}

\begin{proof}
From the compactness of $C$ we conclude that there exists 
at least one pair of points of $\bd (C)$ distant by $\diam (C)$. 
Take any such a pair $f, g$.

By Claim 3 of \cite{L-AEQ} there is a lune $L = H(f) \cap H(g)$ containing $C$. 
Its bounding $(d-1)$-dimensional hemispheres $H(f)/H(g)$ and $H(g)/H(f)$ are centered at $f$ and $g$, respectively.
Moreover, $fg$ is orthogonal to $H(f)/H(g)$ at $f$ and to $H(g)/H(f)$ at $g$.
Of course, $\Delta (L) = |fg|$.
Hence $\Delta (L) = \diam (C)$.

In order to show the first assertion of our claim assume that $\diam (C) = \frac{\pi}{2}$.
If at least one of points $f, g$ is an extreme points of $C$, there is nothing to prove.
Thus further we consider only the situation when both $f$ and $g$ are not extreme.
Since the hemisphere $H(f)/H(g)$ supports $C$, 
we see that an extreme point $e$ of $C$ belongs to $H(f)/H(g)$.
The point $e$ is different from $f$ since $f$ is not an extreme point of $C$. 
Consequently, by the second part of Lemma 3 of \cite{L-AEQ} we see that $|eg| = \frac{\pi}{2}$.
Thus $e, g$ is a promised pair of points.

Let us show the second assertion.
Assume that $\diam (C) < \frac{\pi}{2}$.
Then $\Delta (L) < \frac{\pi}{2}$.
By the first part of Lemma 3 of \cite{L-AEQ} we conclude that every point of $H(f)/H(g)$ different from $f$ is in a distance over $|fg|$ from $g$. 
So since $C$ has an extreme point in $H(f)/H(g)$,
we see that $f$ is an extreme point of $C$.
Analogously, $g$ is an extreme point of $C$.
This confirms the second assertion of our claim.
\end{proof}

By the way, {\it every convex body $C \subset S^2$ of diameter $\frac{\pi}{2}$ contains a pair of extreme points distant by $\frac{\pi}{2}$}.
Here is why. 
Take points $a, b$ from the first assertion of Claim \ref{extreme}, where $b$ extreme. 
If $a$ is not extreme, take this semicircle bounding the lune $L$ from Claim 3 of \cite{L-AEQ}, whose center is $a$.
There are extreme points $a_1, a_2$ of $C$ on this semicircle with $a \in a_1a_2$.
By the second part of Lemma 3 of \cite{L-AEQ} we have $|a_1b| = \frac{\pi}{2}$.

\section{Diameter, width and thickness of reduced spherical bodies}

We say that a convex body $R \subset S^d$ is {\it reduced} if $\Delta (Z) < \Delta (R)$ for each convex body $Z$ being a proper subset of $R$.
Some properties of spherical reduced bodies are given in \cite{L-AEQ}, \cite{LaMu-BULL} \cite{LaMu-AEQ} \cite{LaMu-FAS} and \cite{Mu}. 
This notion is analogous to the notion of a reduced convex body in Euclidean space and finite-dimensional normed space; for instance see the survey articles \cite{LM1} and \cite{LM2}.
Clearly, every spherical body of constant width is a reduced body.
A simple example of a reduced body on $S^2$ is every regular odd-gon.
For more examples see \cite {L-AEQ} and \cite {LaMu-BULL}.   

\begin{thm} \label{Delta-DIAM}
For every reduced body $R\subset S^d$ such that $\Delta(R) \leq \frac{\pi}{2}$ we have $\diam (R) \leq \frac{\pi}{2}$.
Moreover,  if $\Delta(R) < \frac{\pi}{2}$, then $\diam (R) < \frac{\pi}{2}$.
\end{thm}

\begin{proof}
Since $R$ is compact, there are points $p, q \in \bd (R)$ such that $|pq| = \diam (R)$.
Apply Proposition 3.5 of \cite{LaMu-BULL} for $p$.
Thus $R \subset H(p)$. 
Hence $|pq| \leq \frac{\pi}{2}$.
Consequently, $\diam (R) \leq \frac{\pi}{2}$.

Let us show the second part of our theorem.
So now we assume that $\Delta (R) < \frac{\pi}{2}$ and our aim is to show that $\diam (R) < \frac{\pi}{2}$.

Suppose the opposite assertion that $\diam (R) \geq \frac{\pi}{2}$.
Then from the inequality $\diam (R) \leq \frac{\pi}{2}$ 
showed in the proof of the first part of our theorem we obtain $\diam (R) = \frac{\pi}{2}$.

Apply Claim \ref{extreme} for the diameter $\frac{\pi}{2}$.
We conclude that there are two points $p, q \in R$ in the distance $\frac{\pi}{2}$ such that at least one of them, say $p$, is an extreme point of $C$.

Since $p$ is an extreme point of $R$, by Theorem 4 of \cite{L-AEQ} there exists a lune $L_p =I \cap J$ of thickness $\Delta (R)$, containing $R$, where $I, J$ are hemispheres such that $p$ is the center of $I/J$. 
See Figure 1 for $d=3$, where the hemisphere $I$ of $S^3$ is seen as the three-dimensional ball as the view on $S^3 \subset E^4$ ``from outside'' in $E^4$.

\eject
%*********

{\ }

\begin{center}

\includegraphics[width=3.45in]{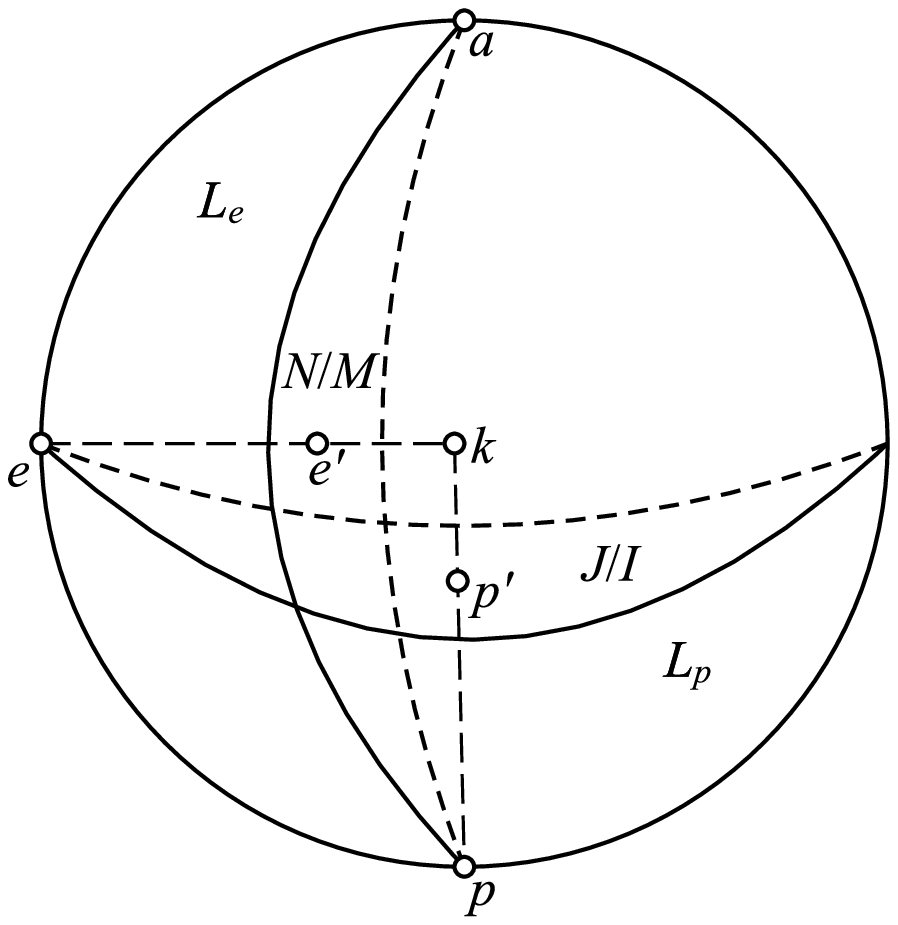} \\ 

{Fig 1. Illustration to the proof of the second statement of Theorem \ref{Delta-DIAM}} 

\end{center}

%**********

By $\Delta (R) < \frac{\pi}{2}$, we have $\Delta (L_p) < \frac{\pi}{2}$.
Thus every point of $L_p$ is in a distance at most $\frac{\pi}{2}$ from $p$.
Hence by $|pq| =  \frac{\pi}{2}$, we conclude that $q$ belongs to ${\rm corn} (L_p)$, which means that the assumption $R \cap {\rm corn} (L_p) \not = \emptyset$ of 
Claim \ref{corners} holds true. 
So by this claim we find an extreme point $e$ of $R$ in ${\rm corn} (L_p)$.
Clearly, $|ep| =  \frac{\pi}{2}$. 

Applying Theorem 4 of \cite{L-AEQ} we find a lune $L_e = M \cap N \supset R$, where $M, N$ are hemispheres, such that $\Delta (L_e) = \Delta (R)$, with $e$ as the center of $M/N$. 

Since $\Delta (R) < \frac{\pi}{2}$, we have $\Delta (L_e) < \frac{\pi}{2}$.
Hence every point of $L_e$ different from its corners is in a distance below $\frac{\pi}{2}$ from $e$.
Thus $p$ must be a corner of $L_e$; just as a point of $R$ in the distance exactly $\frac{\pi}{2}$ from $e$.
We see that the whole arc $pe$ (it is a subset of $R$) is in $M/N$.
So $M/N$ contains also the arc $ea$, where $a$ denotes the antipode of $p$ on $S^d$. 

We see that $M/N$ contains the great semi-circle containing $p$, $e$ and $a$.
Denote the center of $N/M$ by $e'$. 
By Claim 2 of \cite{L-AEQ} we have $e' \in R$.
Clearly, $e' \in \bd (R)$.
From this, since the center $k$ of $I$ does not belong to $N/M$ and since $e'$ is a point of $ek$ different from $p$
we conclude that the center $p'$ of $J/I$ is not in $L_e$. 
Consequently, $L_e$ does not contain the whole $R$, which contradicts the description of $L_e$.

This contradiction shows that our opposite assertion $|pq| \geq \frac{\pi}{2}$ from the third paragraph of this proof is false.
Hence $|pq| < \frac{\pi}{2}$, which means that ${\rm diam} (R) < \frac{\pi}{2}$, which ends the proof of the second thesis of our theorem.
\end{proof}

The special case of the first assertion of Theorem \ref {Delta-DIAM} for $d=2$ is stated in the observation just before Proposition~1 of \cite{LaMu-FAS}.
By the way, the first assertion of Proposition 1 of \cite{LaMu-FAS} is a special case for $d=2$ of the first statement of our Claim \ref {diameter} and the second assertion of Theorem \ref {Delta-DIAM}.
Let us add that our approach is different from proving Proposition~1 of \cite{LaMu-FAS} which applies Theorem 1 of \cite{LaMu-FAS} proved only for $d=2$.

The following theorem generalizes Theorem 4.3 of \cite{LaMu-BULL} from $S^2$ up to $S^d$.
Our proof is ana- 

\noindent
logous (but this time we must apply Proposition 1 of \cite{LaMu-AEQ}).
Here we provide a more detailed consideration supplemented by a figure.

\begin{thm}\label{Delta-CW} 
If a reduced convex body $R \subset S^d$ fulfills  $\Delta(R) \geq \frac{\pi}{2}$, then $R$ is a body of constant width $\Delta (R)$. 
\end{thm}

\begin{proof}
We consider two cases.

Case 1, when $\Delta (R) > \frac{\pi}{2}$.
We apply Proposition 1 of \cite{LaMu-AEQ} that every reduced spherical convex body of thickness over $\frac{\pi}{2}$ is smooth, and next Theorem 5 of \cite{L-AEQ}, that every smooth reduced body is of constant width.

\smallskip
Case 2, when $\Delta (R) = \frac{\pi}{2}$.
Then the assertion of our theorem means that $\width_G (R) = \frac{\pi}{2}$ for every hemisphere $G$ supporting $R$.

In order to confirm this assertion, we intend to get a contradiction under the 
opposite assumption.
Just assume that there iss a hemisphere $K$ supporting $R$ for which 
$\width_K (R)~\not=~\frac{\pi}{2}$. 

By the definition of the thickness, from $\Delta (R) = \frac{\pi}{2}$ we see that $\width_K (R) < \frac{\pi}{2}$ is impossible.
Thus our contrary assertion is nothing else but $\width_K (R) > \frac{\pi}{2}$ (under this form of the contrary assertion, we present Fig. 2 for $d=2$ showing the hemisphere $K$ from the front by the orthogonal outside look at $S^2$).

%**********

\begin{center}

\includegraphics[width=3.1in]{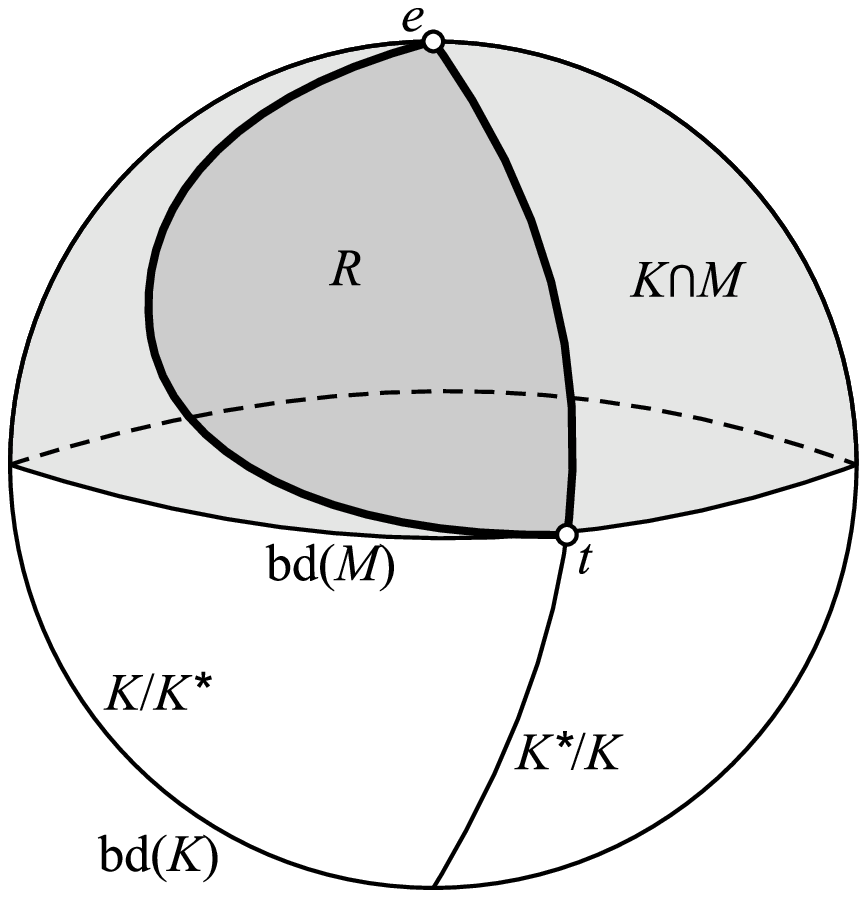} \\ 

\medskip
{Fig. 2. Illustration to the proof of Case 2 of Theorem \ref{Delta-CW}} 

\end{center}

%**********

By Part III of Theorem 1 of \cite{L-AEQ} for our $K$ there exists at least one hemisphere $K^*$ described there.
Take any such $K^*$ and the center $t$ of $K^*/K$ described in this Part III (look also to Corollary 2 there). 
Since $K$ supports $R$, a point $e$ of $R$ belongs to $\bd (K)$. 
Thus $e \in K/K^*$. 
By Proposition 3.5 of \cite{LaMu-BULL} the hemisphere $M$ with center $e$ contains $R$. 
Hence $|te| \le \frac{\pi}{2}$.

By Claim \ref {twolunes} we see that $\bd (K)$ dissects $M$ into two lunes of thickness $\frac{\pi}{2}$.
Since one of them is $K\cap M$, we have $\Delta (K \cap M) = \frac{\pi}{2}$. 
Since $K \cap K^*$ is a narrowest lune over all lunes of the form $K \cap K'$, where the hemisphere $K'$ supports $R$, we have $\Delta(K \cap K^*) \le \Delta (K \cap M)$.
This implies 
$\Delta (K \cap K^*) \le  \frac{\pi}{2}$, which contradicts $\width_K (R) > \frac{\pi}{2}$.  
Therefore our contrary assertion assumed in the second paragraph of Case 2 is false.
Consequently, $\width_G (R) = \frac{\pi}{2}$ for every hemisphere $G$ supporting $R$, which means that
$R$ is of constant width also in Case 2.
\end{proof}

Thanks to Theorem 4 of \cite{LaMu-AEQ}, which says that every spherical body of constant width $w$ has diameter $w$, from Theorem \ref {Delta-CW} we obtain the following proposition.

\begin{pro} \label{diam=Delta} 
For every reduced body $R \subset S^d$ fulfilling $\Delta (R) \geq \frac{\pi}{2}$ we have $\Delta (R) = \diam (R)$. 
\end{pro}

Observe that this proposition is not true without the assumption that the body is reduced.

\begin{cor} 
If $\Delta (R) < \diam (R)$ for a reduced body $R  \subset S^d$, then both these numbers are below $\frac{\pi}{2}$.
Moreover, $R$ is not a body of constant width. 
\end{cor}

\begin{proof}
The inequality $\Delta (R) \geq  \frac{\pi}{2}$ is impossible, since then by Proposition \ref {diam=Delta} we have $\Delta (R) = \diam (R)$, which contradicts the assumption of our statement.
Hence $\Delta (R) <  \frac{\pi}{2}$.
Then by the second part of Theorem \ref {Delta-DIAM} we have $\diam (R) < \frac{\pi}{2}$. 
So both considered numbers are below $\frac{\pi}{2}$.

In order to show the second assertion, assume the opposite that $R$ is of constant width.
Then the maximum and minimum widths of $R$ are equal.
Therefore $\Delta (R) = \diam (R)$. 
This contradicts the assumption of our corollary.
Consequently, $R$ is not of constant width.
\end{proof}

\begin{thm}\label{iff}
Let $R \subset S^d$ be a reduced body. 
Then
\vskip 0.1cm

{\rm (a)} $\Delta(R) =\frac{\pi}{2}$  if and only if $\diam (R) =\frac{\pi}{2}$,

{\rm (b)} $\Delta(R) \geq \frac{\pi}{2}$  if and only if $\diam (R) \geq \frac{\pi}{2}$,

{\rm (c)} $\Delta(R) > \frac{\pi}{2}$  if and only if $\diam (R) > \frac{\pi}{2}$,

{\rm (d)} $\Delta(R) \le \frac{\pi}{2}$  if and only if $\diam (R) \le \frac{\pi}{2}$, 

{\rm (e)} $\Delta(R) < \frac{\pi}{2}$  if and only if $\diam (R) < \frac{\pi}{2}$. 
\end{thm}

\begin{proof}
Let us show the equivalence (a).
If $\Delta(R) =\frac{\pi}{2}$, then by Corollary 1 we have $\diam (R) = \frac{\pi}{2}$.
Now assume that  $\diam (R) =\frac{\pi}{2}$.
Then by the first assertion of Claim \ref{diameter} we get $\Delta (R) \le \frac{\pi}{2}$.
Moreover, by the contrapositive of the second assertion of Theorem \ref {Delta-DIAM} we get $\Delta (R) \geq \frac{\pi}{2}$.
Consequently,  $\Delta (R) = \frac{\pi}{2}$.

In this paragraph we are showing the equivalence (b). 
By Proposition \ref{diam=Delta}, the inequality $\Delta (R) \geq \frac{\pi}{2}$ implies $\diam (R) = \Delta (R)$ and thus $\diam (R) \geq \frac{\pi}{2}$.
The opposite implication is the contrapositive of the second assertion of Theorem \ref{Delta-DIAM}.

From (a) and (b) we get (c). 
It implies (d).  
We obtain (e) as the contrapositive of (b). 
\end{proof}

After Part 4 of \cite{LaMu-AEQ}, we say that a convex body $D \subset S^d$ of diameter $\delta$ is {\it of constant diameter} $\delta$ if for any $p \in \bd (D)$ there exists $p' \in \bd (D)$ such that $|pp'| = \delta$ (more general, this notion makes sense for a closed set $D$ of diameter $\delta$ in a metric space $M$, such that always $x, z \in D$ and $y \in M$ with $|xy| + |yz| = |xz|$ imply $y \in D$, so
for instance when $M$ is a Riemannian manifold). 

We get some spherical bodies of constant diameter on $S^2$ as a particular case of the example given on p. 95 of \cite {LaMu-BULL} by taking there any non-negative $\kappa < \frac{\pi}{2}$ and $\sigma = \frac{\pi}{4} - \frac{\kappa}{2}$. 
The following example presents a wider class of spherical bodies of constant diameter $\frac{\pi}{2}$.

\vskip0.5cm 
\noindent
{\bf Example.}
Take a triangle $v_1v_2v_3 \subset S^2$ of diameter at most $\frac{\pi}{2}$, and put $\kappa_{12} =~|v_1v_2|, \kappa_{23}=|v_2v_3|, \kappa_{31}=|v_3v_1|$,
$\sigma_1 = \frac{\pi}{4} - \frac{\kappa_{12}}{2} + \frac{\kappa_{23}}{2} - \frac{\kappa_{31}}{2}$,
$\sigma_2 = \frac{\pi}{4} - \frac{\kappa_{12}}{2} - \frac{\kappa_{23}}{2} + \frac{\kappa_{31}}{2}$,
$\sigma_3 = \frac{\pi}{4} + \frac{\kappa_{12}}{2} - \frac{\kappa_{23}}{2} - \frac{\kappa_{31}}{2}$.
Here we agree only for triangles with the sum of  lengths of two shortest sides at most the length of the longest side plus $\frac{\pi}{2}$ (equivalently: with $\sigma_1 \geq 0$, $\sigma_2 \geq 0$ and $\sigma_3 \geq 0$).

Prolong the following:
$v_1v_2$ up to $w_{12}w_{21}$ with $v_1 \in w_{12}v_2$,
$v_2v_3$ up to $w_{23}w_{32}$ 
with $v_2 \in w_{23}v_3$, 
$v_3v_1$ up to sides: 
$w_{31}w_{13}$ with $v_3 \in w_{31}v_1$  
(see Fig. 3) such that 
$|v_1w_{12}| = |v_1w_{13}| = \sigma_1$,
$|v_2w_{21}| = |v_2w_{23}| = \sigma_2$,
$|v_3w_{31}| = |v_3w_{32}| = \sigma_3$.
Draw six pieces of circles: 
with center $v_1$ of radius $\sigma_1$ from $w_{12}$ to $w_{13}$ and of radius $\frac{\pi}{2} - \sigma_1$ from $w_{21}$ to $w_{31}$, 
with center $v_2$ of radius $\sigma_2$ from $w_{23}$ to $w_{21}$ and of radius $\frac{\pi}{2} - \sigma_2$ from $w_{32}$ to $w_{12}$, 
with center $v_3$ of radius $\sigma_3$ from $w_{31}$ to $v_{32}$ and of radius $\frac{\pi}{2} - \sigma_3$ from $w_{13}$ to $w_{23}$.
Clearly, the convex hull of these six pieces of circles

{\ }
%**********
\begin{center}

\includegraphics[width=3.1in]{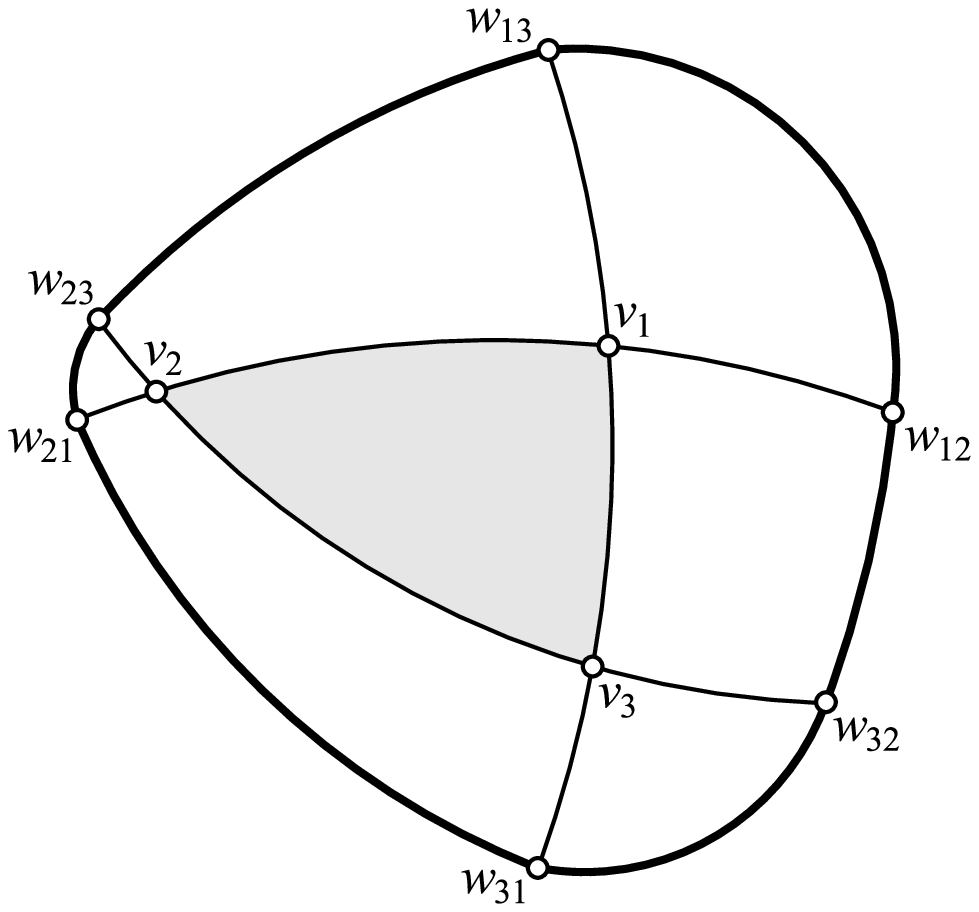} \\ 

\vskip0.1cm
{Fig. 3. A spherical body of constant diameter} 

\end{center}
%********** 

\noindent
is a body of constant diameter~$\frac{\pi}{2}$.

\vskip0.2cm
Generalizing, take a convex odd-gon $v_1 \dots v_n \subset S^2$ of diameter at most $\frac{\pi}{2}$.
 For $i=1, \dots ,n$ put $\kappa_{i \ {i + (n-1)/2}}\ = |v_iv_{i + (n-1)/2}|$ (here and later we mean indices modulo $n$). 
Let $\sigma_i = \frac{\pi}{4} + \Sigma_{i=1}^n s_i\kappa_{i \ {i + (n-1)/2}}$, where $s_i = \frac{1}{2}$ if $|\frac{n+1}{2} -i| \leq \frac {n-3}{4}$ for $n$ of the form $3+4k$ and $|\frac{n}{2} - i| \leq \frac {n-1}{4}$ for $n$ of the form $5+4k$ (where $k=0,1,2, \dots$), and $s_i = -\frac{1}{2}$ in the opposite case.
We agree only for odd-gons with $\sigma_i \geq 0$ for $i=1, \dots , n$.
Prolong each diagonal $v_iv_{i + (n-1)/2}$ up to the arc $w_{i \ i + (n-1)/2} w_{i + (n-1)/2  \ i}$ such that $v_i \in w_{i \ i + (n-1)/2}v_{i + (n-1)/2}$ and $|v_iw_{i \ i + (n-1)/2}| = |v_iw_{i \ i + (n+1)/2}|= \sigma_i$.
For $i= 1, \dots , n$ we draw the piece of the circle with center $v_i$ of radius $\sigma_i$ from $w_{i \ i + (n-1)/2}$ to $w_{i \ i + (n+1)/2}$ and the piece of circle of radius $\frac {\pi}{2} - \sigma_i$ from $w_{i + (n-1)/2  \ i}$ to $w_{i + (n+1)/2  \ i}$. 
The convex hull of the union of our $2n$ pieces of circles is a convex body of constant diameter $\frac{\pi}{2}$.

\vskip0.35cm
Next proposition is applied in the proof of the forthcoming Theorem \ref{wulff}.

\begin{pro} \label{equivalent}
The following conditions are equivalent:

{\rm (1)} \ $C \subset S^d$ is a reduced body with $\Delta(C) = \frac{\pi}{2}$,

{\rm (2)} \ $C \subset S^d$ is a reduced body with $\diam (C) = \frac{\pi}{2}$,

{\rm (3)} \ $C \subset S^d$ is a body of constant width $\frac{\pi}{2}$,

{\rm (4)} \ $C \subset S^d$ is of constant diameter $\frac{\pi}{2}$.
\end{pro}

\begin {proof}
The equivalence of (1) and (2) results from (a) of Theorem \ref{iff}. 
By Theorem \ref {Delta-CW} we conclude that (1) implies (3).
The opposite implication is obvious since every body of constant width is a reduced body.
The equivalence of (3) and (4) follows by the fact (being a particular case of Theorem 5 of \cite{LaMu-AEQ}) that {\it a convex body $W \subset S^d$ is of constant diameter $\frac{\pi}{2}$}, if and only if $W$ is of constant width $\frac{\pi}{2}$.
\end {proof}

\vskip0.2cm
Remark. In particular, Proposition \ref{equivalent} concerns any reduced polygon $V$ of thickness $\frac{\pi}{2}$.
So Theorem 3.2 of \cite{L-COLL} matters for $\Delta (V) = \frac{\pi}{2}$ (instead of $\Delta (V) < \frac{\pi}{2}$ as in \cite{L-COLL}).
The proof is analogous, but let us explain why in lines 16--17 of \cite{L-COLL} the lune $L$ with centers $v_i$ and $t_i$ of its bounding semicircles strictly supports $V$ at $v_i$ (we keep here the notation of \cite{L-COLL} also for $\Delta(V) = \frac{\pi}{2}$). 

Let $p_{i-1}$ (resp. $p_{i+1}$) denote the projection of $v_{i+(n-1)/2}$ (resp. $v_{i+(n+1)/2}$) on $v_iv_{i-1}$ (resp. on $v_iv_{i+1}$). 
Prolong the arc $v_{i+(n-1)/2}p_{i-1}$ (resp. $v_{i+(n+1)/2}p_{i+1}$ ) up to the arc $v_{i+(n-1)/2}r_{i-1}$ (resp. $v_{i+(n+1)/2}r_{i+1}$), where $r_{i-1}$ (resp. $r_{i+1}$) belongs to the semicircle through $v_i$ bounding $L$.
Denote by $z$ the intersection of arcs $v_{i+(n-1)/2}r_{i-1}$ and $v_it_i$. 
We have $|v_{i+(n-1)/2}r_{i-1}| = |v_{i+(n-1)/2}| + |zr_{i-1}|  > |t_iz| + |zv_i| = \frac{\pi}{2}$.  
Analogously, $|v_{i+(n+1)/2}r_{i+1}| > \frac{\pi}{2}$. 
Since $|v_{i+(n-1)/2}p_{i-1}| = \frac{\pi}{2}$ (resp. $|v_{i+1}p_{i+1}| = \frac{\pi}{2}$), we obtain that $v_{i-1}$  (resp. $v_{i+1}$) is in the interior of $L$. 
Hence $L$ strictly supports $V$ at $v_i$.  
So really the whole Theorem 3.2 of \cite{L-COLL} is true also for $\Delta(V) = \frac{\pi}{2}$.

Consequently, also Corollaries 3.6, 3.7, 3.8 and 3.10 of \cite{L-COLL} hold if $\Delta(V) = \frac{\pi}{2}$.

\vskip0.5cm
By part (b) of Theorem \ref{iff}, we conclude the following variant of Theorem \ref{Delta-CW}.

\begin{cor}
If a reduced convex body $R \subset S^d$ fulfills $\, \diam(R) \geq \frac{\pi}{2}$, then $R$ is a body of constant width $w$ equal to $\diam (R)$. 
It is also a body of constant diameter $w$.
\end{cor}

\section{An application for recognizing if a Wulff shape is self-dual}
 
Wulff \cite{Wu} defined a geometric model of a crystal equilibrium, later named {\it Wulff shape}.
The literature concerning this and related subjects is very comprehensive. 
For instance, see the monograph \cite {PV} and the articles 
 \cite {HN1}, \cite{HN3} and 
\cite {TCH}. 

For any continuous function $\gamma : S^d \to  \mathbb{R}_+$, where $\mathbb{R}_+$ denotes the set of positive reals, and $\theta \in S^d$, by $\Gamma_{\gamma, \theta}$ we mean the set of $x \in E^{d+1}$ such that $x \cdot \theta \leq \gamma (\theta)$.
Here the dot means the scalar product of vectors.
The {\it Wulff shape associated with} $\gamma$ is the set ${\cal W}_\gamma = \cap_{\theta \in S^d} \Gamma_{\gamma, \theta}$.
The subject is so well known that we omit here details.

On the other hand, for every convex body $W \subset E^{d+1}$ containing the origin of $E^{d+1}$ in the interior, there exists a continuous function $\gamma : S^d \to \mathbb{R}_+$ such that $W= {\cal W}_\gamma$ (see \cite{Ta}).
Take into account the unique point $(\theta, w(\theta))$ of the intersection of $\bd ({\cal W}_{\gamma})$ with the half-line consisting of points $(\theta , r)$, where $r \in \mathbb{R}_+$.

For a given Wulff shape ${\cal W}_\gamma$ in $E^{d+1}$, Han and Nishimura consider the {\it dual Wulff shape} ${\cal W}_{\bar \gamma}$, where $\bar \gamma (\theta) = 1/w(-\theta)$.
It is denoted by ${\cal DW}_\gamma$.
Next they consider the {\it self-dual Wulff shape} as a Wulff shape ${\cal W}_\gamma$ fulfilling ${\cal W}_\gamma = {\cal D}W_\gamma$.
These notions and their properties and applications are considered
in a number of articles (for instance see \cite{HN1} and \cite{HN3}). 

Han and Nishimura (we follow their notation from \cite{HN1} and \cite{HN3}) apply the classical notion of the central projection $\alpha_N$ from the open hemisphere $H(N)$ centered at a point $N \in S^d$ into the hyperplane $P(N) \subset E^{d+1}$ supporting $S^d$ at $N$. 
This hyperplane may be treated as the $E^d$ with the origin $N$.
The image of a Wulff shape $W_\gamma$ on $P(N)$ 
under the inverse projection $\alpha^{-1}_N$ onto $H(N)$ is called the {\it spherical convex body induced by} $W_\gamma$.
Han and Nishimura \cite{HN1} prove that a Wulff shape $W_\gamma \subset E^d$ is self-dual if and only if the spherical convex body induced by $W_\gamma$ is a spherical body of constant width $\frac{\pi}{2}$.
By this result of them, from Proposition \ref{equivalent} we obtain the following theorem.

\begin{thm} \label{wulff}
{\it Each of the following conditions is equivalent to the statement that the Wulff shape $W_\gamma$ is self-dual:

\noindent
- the spherical convex body induced by $W_\gamma$ is of constant width $\frac{\pi}{2}$,

\noindent
- the spherical convex body induced by $W_\gamma$ is a reduced body of thickness $\frac{\pi}{2}$,

\noindent
- the spherical convex body induced by $W_\gamma$ is a reduced body of diameter $\frac{\pi}{2}$,

\noindent
- the spherical convex body induced by $W_\gamma$ is a body of constant diameter $\frac{\pi}{2}$.}
\end{thm}

Let us add that the equivalence to the third condition  gives the positive answer to the question by Han and Nishimura put at the end of \cite{HN2}.


\baselineskip 10pt 

\begin{thebibliography}{20}

\bibitem{HN1}  
H. Han and T. Nishimura, Self-dual Wulff shapes and spherical convex bodies of constant width $\pi/2$, {\it J. Math. Soc. Japan} {\bf 69} (2017),  1475--1484. 

\bibitem{HN2} 
H. Han and T. Nishimura, Wulff shapes and their duals - RIMS, Kyoto University, http://www.kurims.kyoto-u.ac.jp/$\sim$kyodo/kokyuroku/contents/pdf/2049-04.pdf$\sim$

\bibitem{HN3}
H. Han and T. Nishimura,
Spherical method for studying Wulff shapes and related topics. Singularities in Generic geometry, pp. 1--53 in
Adv. Stud. Pure Math. {\bf 78} (2018) Math. Soc. Japan, Tokyo. 
 
\bibitem{L-AEQ}   
M. Lassak, Width of spherical convex bodies, {\it Aequationes Math.} {\bf 89} (2015), 555--567. 

\bibitem{L-COLL}   
M. Lassak, Reduced spherical polygons, {\it Colloqium Mat.} {\bf 138} (2015), 205--216. 

\bibitem{LM1}   
M. Lassak and H. Martini,
Reduced convex bodies in Euclidean space -- a survey, {\it Expositiones Math.} {\bf 29} (2011), 204--219. 

\bibitem{LM2} 
M. Lassak and H. Martini, Reduced convex bodies in finite-dimensional normed spaces –- a survey, Results Math. {\bf 66} (2014),  
405-426.

\bibitem{LaMu-BULL}   
M. Lassak and M. Musielak, Reduced spherical convex bodies, {\it Bull. Pol. Ac. Math.} {\bf 66} (2018), 87--97. 

\bibitem{LaMu-AEQ}  
M. Lassak and M. Musielak, Spherical bodies of constant width, {\it Aequationes Math.} {\bf 92} (2018), 627--640.

\bibitem{LaMu-FAS}  
M. Lassak and M. Musielak,  Diameter of reduced spherical bodies, {\it Fasciculi Math.}, {\bf 61} (2018), 103--108. 

\bibitem{MMO} 
H. Martini, L. Montejano and D. Oliveros, Bodies of constant width. An introduction to convex geometry with applications, Birkh\"auser, Cham, 2019. 

\bibitem{MS}
H. Martini and K. J. Swanepoel, The geometry of Minkowski spaces —- a survey. Part II. {\it Expo. Math.} {\bf 22} (2004), 93–-144. 

\bibitem{Mu}   
M. Musielak, Covering a reduced spherical body by a disk, {\it Ukr. Math. J.}, to appear (see also arXiv:1806.04246).

\bibitem{Pa}
A. Papadopoulos, On the works of Euler and his followers on spherical geometry, {\it Ganita Bharati} {\bf 36} (2014), no.1, 53--108.

\bibitem{PV}  
A. Pimpinelli and J. Vilain, Physics of Crystal Growth, Monographs and Texts in Statistical Physics, Cambridge University Press, Cambridge, New York, 1998.

\bibitem{Ta} 
J. E. Taylor, Cristalline variational problems, {\it Bull. Amer. Math. Soc.} {\bf 84} (1978), 568--588.

\bibitem{TCH} 
J. E. Taylor, J. E. Cahn and C. A. Handwerker, Geometric models of crystal growth, {\it Acta Metallurgica et Materialia} {\bf 40} (1992), 1443--1474.

\bibitem{Wu}
G. Wulff, Zur Frage der Geschwindigkeit des Wachstrums und der Aufl\"osung der Krystallfl\"achen. In: Zeitschrift f\"ur Krystallographie und Mineralogie {\bf 34}, (1901), 449-530.

\bibitem {VB}  
G. Van Brummelen, Heavenly mathematics. The forgotten art of spherical trigonometry. Princeton University Press, Princeton, 2013. 

\end{thebibliography}
\end{document}